\documentclass[12pt,a4paper]{amsart}
\usepackage[english]{babel}
\usepackage{amsmath,amssymb,amsthm,latexsym,MnSymbol}
\usepackage[utf8]{inputenc}
\usepackage{tikz-cd}
\usepackage[colorlinks=true]{hyperref}
\usepackage[a4paper, portrait]{geometry}
\usepackage{todonotes}
\usepackage{graphicx}

\newcommand{\Z}{\mathbf{Z}}
\newcommand{\Q}{\mathbf{Q}}

\newcommand{\R}{\mathbf{R}}
\newcommand{\C}{\mathbf{C}}

\newcommand{\p}{\mathbf{P}}

\newcommand{\bolda}{{\boldsymbol a}}
\newcommand{\boldb}{{\boldsymbol b}}
\newcommand{\boldc}{{\boldsymbol c}}
\newcommand{\boldd}{{\boldsymbol d}}

\newcommand{\boldx}{{\boldsymbol x}}
\newcommand{\boldy}{{\boldsymbol y}}

\DeclareMathOperator{\dv}{div}
\DeclareMathOperator{\Aut}{Aut}

\DeclareMathOperator{\im}{Im}

\DeclareMathOperator{\SL}{SL}

\newtheorem{thm}{Theorem}
\newtheorem*{thm*}{Theorem}
\newtheorem{lem}[thm]{Lemma}
\newtheorem{pro}[thm]{Proposition}

\newtheorem*{cor*}{Corollary}
\newtheorem{conjecture}[thm]{Conjecture}
\theoremstyle{definition}
\newtheorem{definition}[thm]{Definition}

\theoremstyle{remark}

\begin{document}

\author[F. Brunault]{Fran\c{c}ois Brunault}

\date{\today}

\address{\'ENS Lyon, Unit\'e de math\'ematiques pures et appliqu\'ees,
46 all\'ee d'Italie, 69007 Lyon, France}

\email{francois.brunault@ens-lyon.fr}
\urladdr{https://perso.ens-lyon.fr/francois.brunault}

\title[On the Mahler measure of $(1+x)(1+y)+z$]{On the Mahler measure of $(1+x)(1+y)+z$}


\thanks{The author was supported by the research project ``Motivic homotopy, quadratic invariants and diagonal classes'' (ANR-21-CE40-0015) operated by the French National Research Agency (ANR)}

\begin{abstract}
We prove a conjecture of Boyd and Rodriguez Villegas relating the Mahler measure of the polynomial $(1+x)(1+y)+z$ and the value at $s=3$ of the $L$-function of an elliptic curve of conductor $15$. The proof makes use of the computation by Zudilin and the author of the regulator of certain $K_4$ classes on modular curves.
\end{abstract}

\maketitle

\section{Introduction}

The (logarithmic) Mahler measure of a complex Laurent polynomial $P(x_1,\ldots,x_n)$ is defined as the average of $\log |P|$ over the torus $T^n : |x_1| = \ldots = |x_n| = 1$,
\begin{equation*}
m(P) = \frac{1}{(2\pi i)^n} \int_{T^n} \log |P| \; \frac{dx_1}{x_1} \ldots \frac{dx_n}{x_n}.
\end{equation*}
For a one-variable polynomial $\displaystyle P(x) = c \, \prod_{i=1}^d (x-\alpha_i)$, Jensen's formula gives
\begin{equation*}
m(P) = \log |c| + \sum_{\substack{i = 1 \\ |\alpha_i| \geq 1}}^d \log |\alpha_i|.
\end{equation*}
Originally, the Mahler measure was introduced as a height function for multivariate polynomials, in relation with transcendental number theory. It was later realised that the Mahler measure appears naturally in other contexts. For example, Smyth discovered in 1981 the following formulas
\begin{equation} \label{eq Smyth}
m(1+x+y) = \frac{3 \sqrt{3}}{4\pi} L(\chi_{-3}, 2), \qquad \qquad m(1+x+y+z) = \frac{7}{2\pi^2} \zeta(3),
\end{equation}
where $\chi_{-3}(n) = (\frac{-3}{n})$ is the Dirichlet character modulo $3$. The Mahler measure of integer polynomials turns out to have deep links with special values of $L$-functions. We mention here some aspects of this connection, referring to \cite{BLRVD03, Boy05, BL13, BZ20} for more complete surveys.

A combination of experiments and theoretical insights led Boyd and Deninger to conjecture the identity
\begin{equation} \label{eq BD}
m\Bigl(x+\frac{1}{x}+y+\frac{1}{y}+1\Bigr) = L'(E,0),
\end{equation}
where $E : x+\frac{1}{x}+y+\frac{1}{y}+1 = 0$ is an elliptic curve of conductor $15$. This was proved some 15 years later by Rogers and Zudilin \cite{RZ14}.

The identity \eqref{eq BD} can be conceptually explained using Beilinson's theory of regulators, and Deninger gave in \cite{Den97} a general framework to relate Mahler measures and cohomology. More precisely, let $P(x_1,\ldots,x_n)$ be a complex Laurent polynomial, which we assume to be monic in $x_n$. Applying Jensen's formula with respect to $x_n$, we may write $m(P)$ as an integral
\begin{equation} \label{eq int eta}
\int_\Gamma \eta(x_1,\ldots,x_n),
\end{equation}
where $\eta$ is a differential $(n-1)$-form on the zero locus $V_P$ of $P$ in $(\C^\times)^n$, and $\Gamma$ is the $(n-1)$-dimensional \emph{Deninger chain},
\begin{equation*}
\Gamma = \bigl\{ (x_1, \ldots, x_n) \in V_P : \; |x_1| = \cdots = |x_{n-1}| = 1, |x_n| \geq 1 \bigr\},
\end{equation*}
endowed with the orientation coming from that of $T^{n-1}$. We make here all necessary assumptions for this integral to make sense \cite[Assumptions 3.2]{Den97}, in particular $\Gamma$ must avoid the singular points of $V_P$.

Assume now that $\Gamma$ is closed. Then \eqref{eq int eta} can be given a cohomological interpretation, since the class of $\eta$ in de Rham cohomology is the image under the Beilinson regulator map of the cup-product $\{x_1,\ldots,x_n\}$ in the motivic cohomology group $H^n_{\mathcal{M}}(V_P, \Q(n))$. This situation is favourable and under certain conditions, the Beilinson conjectures predict a link between $m(P)$ and some $L$-value associated to $V_P$. The identity \eqref{eq BD} is an example of this phenomenon (in reality, in this case the path $\Gamma$ is not closed, but symmetries can be used to ``close the path'').

A more mysterious situation is when the form $\eta$ is exact, in which case we say $P$ is \emph{exact}. Stokes's formula reduces the Mahler measure $m(P)$ to an $(n-2)$-dimensional integral over the boundary $\partial \Gamma$, but Deninger's theory does not provide an intrinsic cohomological interpretation of this integral. Maillot suggested in 2003 that, in the exact case, $m(P)$ should be related to the cohomology of the variety
\begin{equation*}
W_P : P(x_1,\ldots,x_n) = \overline{P}\Bigl(\frac{1}{x_1},\ldots,\frac{1}{x_n}\Bigr) = 0.
\end{equation*}
What makes it plausible is that $\partial \Gamma$ is contained in $W_P$, because $V_P \cap T^n = W_P \cap T^n$. The relevant motivic cohomology group is now $H^{n-1}_{\mathcal{M}}(W_P, \Q(n))$, which is harder to deal with, as we cannot use cup-products. The identities \eqref{eq Smyth} are of this type. For example, the polynomial $1+x+y$ leads to the algebraic $K$-group $K_3(\Q(\sqrt{-3}))$, which is known to have rank $1$ by Borel's theorem. In general, motivic cohomology $H^i_{\mathcal{M}}(\cdot, \Q(n))$ with $i \neq n$ makes it more difficult to handle the Mahler measure.

Following Maillot's insight, Boyd and Rodriguez Villegas discovered in 2003 several identities involving $3$-variable exact polynomials \cite{BLRVD03, Boy05, Boy06a}. One example is:

\begin{conjecture}[Boyd and Rodriguez Villegas \cite{Boy05}] \label{conj BRV}
We have the equality
\begin{equation} \label{eq conj BRV}
m((1+x)(1+y)+z) \stackrel{?}{=} -2 L'(E,-1),
\end{equation}
where $E : (1+x)(1+y)(1+\frac{1}{x})(1+\frac{1}{y})=1$ is an elliptic curve of conductor $15$.
\end{conjecture}

Here $E$ arises as the Maillot variety of $P = (1+x)(1+y)+z$. The first result towards Conjecture \ref{conj BRV} was obtained by Lal\'\i n \cite{Lal15}, who related the Mahler measure of $P$ to the regulator of a cocycle in the Goncharov complex $\Gamma(E,3)$ (see Section \ref{sec weight 3} for the definition of this complex). Let us write $\gamma_E = \partial \Gamma$ for the boundary of Deninger's chain $\Gamma$; this is a closed path in $E$.

\begin{thm}[Lal\'\i n] \label{thm Lalin}
We have $m(P) = \frac{1}{4\pi^2} \int_{\gamma_E} r_3(2)(\xi_E)$, where $\xi_E$ is the class of the cocycle $\{-x\}_2 \otimes y - \{-y\}_2 \otimes x$ in $\Gamma(E, 3)$, and $r_3(2)$ is the Goncharov regulator map.
\end{thm}

In essence, Lal\'\i n's theorem reduces Conjecture \ref{conj BRV} to the Beilinson conjecture for $L'(E,-1)$. In this article, we compute the above Goncharov regulator, leading to the following theorem.

\begin{thm} \label{thm main}
The Boyd and Rodriguez Villegas conjecture \eqref{eq conj BRV} is true.
\end{thm}

Another fascinating conjecture by Rodriguez Villegas concerns the Mahler measure of the polynomials $1+x_1+\ldots+x_n$ for $n = 4$ and $n = 5$. These polynomials are also exact and their Mahler measures are expected to involve $L$-values of cusp forms of weight $3$ and $4$, respectively \cite[Section 6.2]{BZ20}. Partial results have been obtained by Shinder and Vlasenko \cite{SV14}. Here is a similar identity that we found recently:
\begin{equation*}
m((1+x)(1+y)(1+z)+t) \stackrel{?}{=} -6 L'(f_7, -1) - \frac{48}{7} \zeta'(-2),
\end{equation*}
where $f_7(\tau) = \eta(\tau)^3 \eta(7\tau)^3$ is the unique CM newform of weight $3$ and level $7$.

The main ingredient in the proof of Theorem \ref{thm main} is the computation by Zudilin and the author \cite{BZ23} of the Goncharov regulator of explicit classes $\xi_1(a,b)$ in the motivic cohomology of the modular curve $Y_1(N)$, which were introduced in \cite{Bru20b}. A key fact here is that $E$ is isomorphic to the modular curve $X_1(15)$, something we make precise in Section \ref{sec parametrisation}. In Section \ref{sec weight 3}, we recall Goncharov's theory of polylogarithmic complexes in weight $2$ and $3$ and, for modular curves, we define subcomplexes built out of modular units. These complexes are amenable to computation, and we partly implemented the weight $3$ complex in PARI/GP \cite{PARI}; the scripts are available at \cite{Bru23}. In Sections \ref{sec class} and \ref{sec path}, we express Lal\'\i n's class $\xi_E$ and the path $\gamma_E$ in purely modular terms. The final computation is performed in Section \ref{sec final}, using the results of \cite{BZ23}. In the appendix, we give tables of (conjectural) identities relating $3$-variable Mahler measures and $L(E,3)$ for a number of elliptic curves $E$ over $\Q$.

\textbf{Acknowledgements.} I am grateful to Matilde Lal\'\i n, Riccardo Pengo, Wadim Zudilin and the International Groupe de travail on differential equations in Paris for exchanges which have been helpful in several parts of this paper. I would also like to thank Berend Ringeling for checking numerically several Mahler measure identities from the appendix.

\section{The modular parametrisation} \label{sec parametrisation}

Consider the polynomial $P(x,y,z) = (1+x)(1+y)+z$. We keep the same notations as in the introduction, so that the Maillot variety $W_P$ in $(\C^\times)^3$ is defined as
\begin{equation*}
W_P : \begin{cases} (1+x)(1+y) + z = 0, \\
(1+\frac{1}{x})(1+\frac{1}{y}) + \frac{1}{z} = 0.
\end{cases}
\end{equation*}
Eliminating $z$, we see that $W_P$ is isomorphic to the smooth curve in $(\C^\times)^2$ given by
\begin{equation} \label{eq C}
C : (1+x)^2 (1+y)^2 = xy.
\end{equation}
Let $E$ denote the closure of $C$ in $\p^1(\C) \times \p^1(\C)$. We view $E$ as a smooth projective curve defined over $\Q$. It turns out that $E$ is isomorphic to an elliptic curve of conductor $15$ \cite[(4.2)]{Lal15}. The PARI/GP commands

\;

\begin{verbatim}
E = ellfromeqn((1+x)^2*(1+y)^2-x*y)
ellidentify(ellinit(E))
\end{verbatim}
confirm that $E$ is isomorphic to the elliptic curve with Cremona label $15a8$. On the other hand, we know that the modular curve $X_1(15)$ is isomorphic to $15a8$, since they are both elliptic curves of conductor $15$, and the period lattice of $X_1(15)$ can be computed using modular symbols, agreeing with that of $15a8$. Note that Stevens's conjecture \cite[Conjecture II]{Ste89} is known in this case by \cite[Section 7]{Ste89}.

In Proposition \ref{pro parametrisation} below we give an explicit isomorphism $X_1(15) \cong E$ (note that the proof does not rely on floating point computations). An important feature of this parametrisation is that the functions $-x$ and $-y$ correspond to modular units on $X_1(15)$. This is crucially used in Section~\ref{sec class} to relate Lal\'\i n's class $\xi_E$ and the modular classes $\xi_1(a,b)$ from \cite[Section 6]{Bru20b}. Even more, we need the functions $-x$ and $-y$ to be of the form $u_1(a,b,c,d)$, a class of modular units introduced in \cite{Bru20b} and whose definition we now recall.

Let $N \geq 1$ be an integer. For any $\bolda = (a_1, a_2) \in (\Z/N\Z)^2/\pm 1$, $\bolda \neq (0,0)$, we define
\begin{equation*}
\wp_\bolda(\tau) = \wp\Bigl(\tau ; \frac{a_1 \tau + a_2}{N}\Bigr) \qquad (\tau \in \C, \, \im(\tau)>0),
\end{equation*}
where $\wp(\tau ; z)$ is the Weierstra\ss{} function. The function $\wp_\bolda$ is a modular form of weight $2$ on the principal congruence group $\Gamma(N)$. For any distinct elements $\bolda, \boldb, \boldc, \boldd$ of $(\Z/N\Z)^2/\pm 1$, we then define $u(\bolda, \boldb, \boldc, \boldd)$ as the cross-ratio $[\wp_\bolda, \wp_\boldb, \wp_\boldc, \wp_\boldd]$. This is a modular unit on $\Gamma(N)$. For distinct elements $a,b,c,d$ of $(\Z/N\Z)/\pm 1$, we use the shortcut $u_1(a,b,c,d) = u((0,a),(0,b),(0,c),(0,d))$, which is a modular unit on $\Gamma_1(N)$. These units are defined over $\Q$. The properties of $u(\bolda, \boldb, \boldc, \boldd)$ needed in this article can be found in \cite[Section 3]{Bru20b}.

In the following proposition, we take $N=15$.

\begin{pro} \label{pro parametrisation}
The curve $E$ is parametrised by the following modular units on $\Gamma_1(15)$:
\begin{equation} \label{eq xtau ytau}
x(\tau) = - u_1(1,2,3,7)(\tau), \qquad y(\tau) = - u_1(2,4,6,1)(\tau).
\end{equation}
Moreover, the map $\tau \mapsto (x(\tau),y(\tau))$ induces an isomorphism $\varphi : X_1(15) \xrightarrow{\cong} E$ defined over $\Q$.
\end{pro}

\begin{proof}
Let us show that $u = - u_1(1,2,3,7)$ and $v = - u_1(2,4,6,1)$ satisfy $(1+u)^2 (1+v)^2 = uv$. For this  we may replace $u$ and $v$ by their transforms under the Atkin-Lehner involution $W_{15} : \tau \mapsto -1/15\tau$ on $X_1(15)$, as this does not affect the equation. The units $\tilde{u} = u \circ W_{15}$ and $\tilde{v} = v \circ W_{15}$ can be expressed in terms of Siegel units of level $15$ using \cite[eq.~(6)]{Bru20b}:
\begin{equation} \label{eq ut vt}
\begin{split}
\tilde{u} & = - \frac{\tilde{g}_2 \tilde{g}_4}{\tilde{g}_1 \tilde{g}_7} = -1 - q + q^4 + q^5 - q^7 + O(q^8), \\
\tilde{v} & = - \frac{\tilde{g}_4 \tilde{g}_7}{\tilde{g}_1 \tilde{g}_2} = -q^{-2} - q^{-1} - 2 - 2q - 2q^2 - 2q^3 - 2q^4 + O(q^5)
\end{split}
\end{equation}
where, for $a \in \Z/N\Z$, $a \neq 0$,
\begin{equation*}
\tilde{g}_a(\tau) = q^{N B_2(\hat{a}/N)/2} \prod_{\substack{n \geq 1 \\ n \equiv a \bmod{N}}} (1-q^n) \prod_{\substack{n \geq 1 \\ n \equiv -a \bmod{N}}} (1-q^n) \qquad (q = e^{2\pi i \tau}).
\end{equation*}
Here $B_2(t)=t^2-t+\frac16$ is the Bernoulli polynomial and $\hat{a}$ is the lift of $a$ in $\{1,\ldots,N-1\}$.

We are now going to compute the divisors of $\tilde{u}$ and $\tilde{v}$. To this end, we recall the description of the cusps of the modular curve $X_1(N)$. There is a bijection \cite[Example 9.1.3]{DI95}
\begin{equation*}
\{\textrm{cusps of } X_1(N)(\C)\} \xrightarrow{\cong} \{(c,d) : c \in \Z/N\Z, \, d \in (\Z/(c,N)\Z)^\times\}/ \pm 1
\end{equation*}
which associates to a cusp $\gamma \infty$ with $\gamma \in \SL_2(\Z)$, the class of the bottom row $(c,d)$ of $\gamma$. Moreover, by \cite[Section 9.3, p.~79]{DI95}, the Galois action on the cusps is described as follows: for $\sigma \in \Aut(\C)$, we have $\sigma \cdot (c,d) = (c, \chi(\sigma) d)$, where $\chi(\sigma) \in (\Z/N\Z)^\times$ is characterised by $\sigma(e^{2\pi i/N}) = e^{2\pi i \chi(\sigma)/N}$. As a consequence, a complete set of representatives of the Galois orbits is provided by the cusps $\frac{1}{k} = (\begin{smallmatrix} 1 & 0 \\ k & 1 \end{smallmatrix}) \infty$ with $0 \leq k \leq \lfloor \frac{N}{2} \rfloor$.

Now we can compute the divisor of $u_1(a,b,c,d)$ for distinct $a,b,c,d \in (\Z/N\Z)/\pm 1$ as follows. Since this unit is defined over $\Q$, it suffices to determine its order of vanishing at the cusps $1/k$ just described. By \cite[Proposition 3.6]{Bru20b}, we have
\begin{equation*}
u_1(a,b,c,d) | \begin{pmatrix} 1 & 0 \\ k & 1 \end{pmatrix} = u((ka,a),(kb,b),(kc,c),(kd,d)).
\end{equation*}
The order of vanishing of this unit at $\infty$ is deduced from the expression of $u(\bolda, \boldb, \boldc, \boldd)$ in terms of Siegel units \cite[Proposition 3.7]{Bru20b}, taking into account that it should be computed with respect to the uniformising parameter $q^{(k,N)/N}$. Applying this in our situation, we obtain
\begin{equation*}
\dv(u) = -2 [1/2] + 2 [1/7], \qquad \qquad \dv(v) = -2[0] + 2[1/4].
\end{equation*}
These cusps are rational and we see that $F = (1+u)^2 (1+v)^2 - uv$ has poles of order at most $4$ at $0$ and $1/2$, and is regular elsewhere. Moreover, we compute from \eqref{eq ut vt} that $F(-1/15\tau) = O(q^5)$ when $\mathrm{Im}(\tau) \to +\infty$. Therefore $F$ vanishes at order $\geq 5$ at $0$, and consequently $F=0$.

It remains to show that $\varphi : X_1(15) \to E$ is an isomorphism. Since $X_1(15)$ and $E$ are smooth, it suffices to check that $\varphi$ is a birational map. We know that $u$ has degree $2$ as a function on $X_1(15)$, while $x$ has degree $2$ as a function on $E$. It follows that $\varphi$ is birational.
\end{proof}

\section{The weight $3$ complex of the modular curve $Y_1(N)$} \label{sec weight 3}

Goncharov has defined in \cite{Gon95} polylogarithmic complexes which are expected to compute the motivic cohomology of arbitrary fields. We define in this section a modular complex $\mathcal{C}_N(3)$, which is a subcomplex of the weight $3$ polylogarithmic complex attached to the modular curve $Y_1(N)$. It is generated (in a suitable sense) by the Siegel units and the modular units $u_1(a,b,c,d)$ from Section \ref{sec parametrisation}. Our construction can be seen as a weight $3$ analogue of the weight~$2$ Euler complex $\mathbb{E}_N^\bullet$ introduced by Goncharov in \cite{Gon08}. We also explain how to manipulate $\mathcal{C}_N(3)$ using PARI/GP. The constructions below work with no more effort for the modular curve $Y(N)$ with full level $N$ structure, using parameters in $(\Z/N\Z)^2$ instead of $\Z/N\Z$. However we have not implemented it, as the case of $Y_1(N)$ suffices for our application.

We briefly recall Goncharov's polylogarithmic complexes in weight $2$ and $3$. Let $F$ be any field. Define $B_2(F)$ to be the quotient of $\Q[F^\times \backslash \{1\}]$ by the subspace generated by the $5$-term relations \cite[Section 1.8]{Gon95}. The group $B_3(F)$ is defined similarly as the quotient of $\Q[F^\times \backslash \{1\}]$ by explicit relations \cite[Section 1.8]{Gon95}, whose definition will not be needed here. For $x \in F^\times \backslash \{1\}$ and $n \in \{2,3\}$, we denote by $\{x\}_n$ the image of the generator $[x]$ in $B_n(F)$. Then the complex $\Gamma(F,2)$, in degrees $1$ and $2$, is defined as
\begin{equation*}
\Gamma(F,2) :
\begin{tikzcd}[row sep = .3em]
B_2(F) \arrow{r} & \Lambda^2 F^\times \otimes \Q \\
 \{x\}_2 \arrow[mapsto]{r} & (1-x) \wedge x,
\end{tikzcd}
\end{equation*}
and the complex $\Gamma(F,3)$, in degrees $1$ to $3$, is defined as
\begin{equation*}
\Gamma(F,3) :
\begin{tikzcd}[row sep = .3em]
B_3(F) \arrow{r} & B_2(F) \otimes F^\times \otimes \Q \arrow{r} & \Lambda^3 F^\times \otimes \Q \\
\{x\}_3 \arrow[mapsto]{r} & \{x\}_2 \otimes x & \\
& \{x\}_2 \otimes y \arrow[mapsto]{r} & (1-x) \wedge x \wedge y.
\end{tikzcd}
\end{equation*}

Goncharov conjectures that $H^i(\Gamma(F,n))$ is isomorphic to $H^i_{\mathcal{M}}(F, \Q(n))$. In the case $F$ is the function field of a smooth curve $Y$ over a field $k$, these complexes are endowed with residue maps $\Gamma(F,n) \to \Gamma(k(x), n-1)[-1]$ for every closed point $x \in Y$. Goncharov then defines the complex $\Gamma(Y,n)$ as the simple of the morphism of complexes $\Gamma(F,n) \to \bigoplus_{x \in Y} \Gamma(k(x), n-1)[-1]$, and he conjectures that $H^i(\Gamma(Y,n))$ is isomorphic to $H^i_{\mathcal{M}}(Y, \Q(n))$ \cite[Section 1.15(b)]{Gon95}.

We will consider these complexes in the case $Y$ is the modular curve $Y_1(N)$, and $F$ is its function field. We will see, in particular, that they have natural subcomplexes built out of modular units.

\begin{definition}
Fix an integer $N \geq 1$. We introduce the following sets of modular units on $Y_1(N)$:
\begin{itemize}
\item $U_1$ consists of the Siegel units $g_{0,a}$, $a \in (\Z/N\Z) \backslash \{0\}$, in $\mathcal{O}(Y_1(N))^\times \otimes \Q$;
\item $U_2$ consists of the modular units $u_1(a,b,c,d)$ in $\mathcal{O}(Y_1(N))^\times$, where $a,b,c,d$ are distinct elements of $(\Z/N\Z)/\pm 1$.
\end{itemize}
Moreover, we associate to them the following spaces:
\begin{itemize}
\item $\langle U_1 \rangle$ is the $\Q$-span of $U_1$ in $F^\times \otimes \Q$;
\item $\langle U_2 \rangle$ is the $\Q$-span of $\{u\}_2$, $u \in U_2$, in $B_2(F)$;
\item $\langle U_2 \rangle_3$ is the $\Q$-span of $\{u\}_3$, $u \in U_2$, in $B_3(F)$.
\end{itemize}
\end{definition}

With these definitions, the weight $2$ modular complex can be defined as
\begin{equation*}
\mathcal{C}_N(2) :
\begin{tikzcd}
\langle U_2 \rangle \arrow{r} & \Lambda^2 \langle U_1 \rangle, \qquad \{u\}_2 \mapsto (1-u) \wedge u.
\end{tikzcd}
\end{equation*}
This complex is well-defined because $U_2$ is contained in $\langle U_1 \rangle$ by \cite[Proposition 3.8]{Bru20b}, and $U_2$ is stable under $u \mapsto 1-u$ from the definition of $u_1(a,b,c,d)$ as a cross-ratio. It would be interesting to compare $\mathcal{C}_N(2)$ with the Euler complex $\mathbb{E}_N^{\bullet}$ defined by Goncharov in \cite[Section 2.5]{Gon08}.

We are now ready to introduce a version of the weight $3$ modular complex.

\begin{definition}
The complex $\mathcal{C}_N(3)$ is the following subcomplex of $\Gamma(F, 3)$ in degrees $1$ to $3$:
\begin{equation*}
\mathcal{C}_N(3) :
\begin{tikzcd}
\langle U_2 \rangle_3 \arrow{r} & \langle U_2 \rangle \otimes \langle U_1 \rangle \arrow{r} & \Lambda^3 \langle U_1 \rangle.
\end{tikzcd}
\end{equation*}
\end{definition}

We warn the reader that the group $\langle U_2 \rangle_3$ in degree $1$ may not be the right one. Indeed, the unit $u_1(a,b,c,d)$ is by definition a cross-ratio, hence is a natural argument for the dilogarithm, but \emph{a priori} not for the trilogarithm. However, the complex $\mathcal{C}_N(3)$ will suffice for our needs.

Since the construction of $\mathcal{C}_N(3)$ involves only modular units, the elements of $\langle U_2 \rangle_3$, $\langle U_2 \rangle \otimes \langle U_1 \rangle$ and $\Lambda^3 \langle U_1 \rangle$ have trivial residues at every point of $Y_1(N)$. In particular, $\mathcal{C}_N(3)$ embeds as a subcomplex of $\Gamma(Y_1(N),3)$, and we have natural maps in cohomology $H^i(\mathcal{C}_N(3)) \to H^i(\Gamma(Y_1(N),3))$ in degree $i \in \{1,2,3\}$. The case of interest to us is $i=2$.

We have implemented part of the complex $\mathcal{C}_N(3)$ in PARI/GP, with the specific aim of comparing cocycles in degree $2$. Firstly, the following lemma gives a natural way to represent modular units in $\langle U_1 \rangle$.

\begin{lem} \label{lem basis U1}
A basis of $\langle U_1 \rangle$ is given by the Siegel units $g_{0,a}$ with $1 \leq a \leq \lfloor N/2 \rfloor$.
\end{lem}

\begin{proof}
We have $g_{0,-a} = g_{0,a}$ in $F^\times \otimes \Q$, and by \cite{Str15}, the units $g_{0,a}$ with $1 \leq a \leq \lfloor N/2 \rfloor$ form a basis of $(\mathcal{O}(Y_1(N))^\times / \Q^\times) \otimes \Q$.
\end{proof}

Each unit in $U_2$ can be written in the basis of Lemma \ref{lem basis U1} using \cite[Proposition 3.8]{Bru20b}. Note that no computation of divisor is needed here, thanks to this choice of basis. We actually need to determine $U_2$ as a set, and so to check whether two given units $u_1(a,b,c,d)$ and $u_1(a',b',c',d')$ are equal. We remark that the leading coefficient of $u_1(a,b,c,d)$ at the cusp $0$ is equal to $1$ by the discussion after \cite[Proposition 3.8]{Bru20b}. Combining this with Lemma \ref{lem basis U1}, we see that two units are equal if and only if their coordinates in the basis of $\langle U_1 \rangle$ are equal.

We now consider the free vector space $\Q[U_2]$, and we quotient it by the following subspaces, encoding the relations between the symbols $\{u_1(a,b,c,d)\}_2$. From the definition of $u_1(a,b,c,d)$ as a cross-ratio, the symmetric group $\mathcal{S}_4$ acts on $U_2$ by permuting the indices, and this action factors through $\mathcal{S}_3$. Moreover, because of the relations $\{1/u\}_2 = \{1-u\}_2 = - \{u\}_2$ in $B_2(F)$ \cite[VI, Lemma 5.4]{Wei13}, we have the antisymmetry property:
\begin{equation} \label{eq antisymmetry}
\{u_1(a_{\sigma(1)}, a_{\sigma(2)}, a_{\sigma(3)}, a_{\sigma(4)})\}_2 = \varepsilon(\sigma) \{u(a_1, a_2, a_3, a_4)\}_2 \qquad (\sigma \in \mathcal{S}_4),
\end{equation}
for all distinct parameters $a_i$ in $(\Z/N\Z)/\pm 1$, where $\varepsilon(\sigma) = \pm 1$ is the signature. It thus suffices to consider those parameters satisfying $0 \leq a < b < c < d \leq \lfloor N/2 \rfloor$. The elements $\{u_1(a,b,c,d)\}_2$ are also subject to the $5$-term relations \cite[Lemma 4.7]{Bru20b}:
\begin{equation} \label{eq 5terms}
\sum_{j \in \Z/5\Z} \{u_1(a_j, a_{j+1}, a_{j+2}, a_{j+3})\}_2 = 0 \qquad \textrm{in } B_2(F),
\end{equation}
for any family $(a_j)_{j \in \Z/5\Z}$ of distinct elements of $(\Z/N\Z)/\pm 1$. We denote by $R_2$ the subspace of $\Q[U_2]$ generated by the antisymmetry relations \eqref{eq antisymmetry} and the $5$-term relations \eqref{eq 5terms}. Finally, we denote by $Q$ the space of degree $2$ coboundaries in the complex $\mathcal{C}_N(3)$, namely, the subspace of $\Q[U_2] \otimes \langle U_1 \rangle$ generated by the symbols $[u] \otimes u$ with $u \in U_2$.

In practice, in order to reduce the size of the objects, we only compute:
\begin{itemize}
\item a set $U'_2$ of representatives of the quotient $U_2 / \mathcal{S}_3$;
\item the subspace $R'_2$ of $\Q[U'_2]$ generated by the $5$-term relations;
\item the subspace $Q'$ of $\Q[U'_2] \otimes \langle U_1 \rangle$ of degree $2$ coboundaries.
\end{itemize}

The corresponding scripts are contained in the file \texttt{K4-modular-complex.gp} from \cite{Bru23}. They can be applied in the following way. Say we have two degree $2$ cocycles $\xi$ and $\xi'$ in $\Gamma(Y_1(N),3)$. Assume that they are both linear combinations of symbols $\{u_1(a,b,c,d)\}_2 \otimes g_{0,x}$. We may then represent $\xi-\xi'$ by an element of $\Q[U_2] \otimes \langle U_1 \rangle$, and we check whether this element belongs to the subspace $R_2 \otimes \langle U_1 \rangle + Q$. If so, then we can deduce that $\xi$ and $\xi'$ are cohomologous, and thus have the same image in $K_4^{(3)}(Y_1(N))$ under De Jeu's map \cite[Theorems 5.3 and 5.4]{Bru20b}. If $\xi-\xi'$ does not belong to the subspace, we cannot conclude anything, as $R_2$ and $Q$ may not contain all the relations in the respective groups.

The linear system involved in the above computation has size $O(N^5) \times O(N^6)$. Experimentally, we have found that the cardinality of $U_2$ for $N=p$ prime is $(p^2-1)(p^2-25)/192$, which is smaller by a factor of about $3$ than what we could expect, namely $6 \binom{(p+1)/2}{4}$. Furthermore, it seems that the dimension of $\Q[U_2]/R_2$ is equal to $(p-1)(p-5)/12$, which is also the number of triples $(a,b,c)$ with $0 < a < b < c < p$ and $a + b + c \equiv 0 \mod{p}$, where $2b < p$; see \cite[Sequence \href{https://oeis.org/A242090}{A242090}]{OEIS}. If true, there should be a way to bypass the step of quotienting by $R_2$. This would result in a much smaller linear system for the comparison of cocycles.

\section{The Lal\'\i n class} \label{sec class}

Recall that Lal\'\i n's theorem (Theorem \ref{thm Lalin}) expresses the Mahler measure of $(1+x)(1+y)+z$ as the regulator integral of the following cocycle in the weight $3$ Goncharov complex of $E$:
\begin{equation*}
\xi_E = \{-x\}_2 \otimes y - \{-y\}_2 \otimes x.
\end{equation*}
Our aim in this section is to relate $\xi_E$ to the classes $\xi_1(a,b)$ on $X_1(15)$, which were introduced in \cite[Section 6]{Bru20b}. This is a purely algebraic computation making use of our implementation of the weight $3$ complex of $X_1(15)$ explained in Section \ref{sec weight 3}.

We first pull back $\xi_E$ to the modular curve $X_1(15)$ using the modular parametrisation $\varphi$. Using Proposition \ref{pro parametrisation} and its proof, we have in the degree $2$ cohomology of $\Gamma(Y_1(15),3)$
\begin{equation} \label{eq xi_15}
\varphi^* \xi_E = \{u_1(1,2,3,7)\}_2 \otimes \Bigl(\frac{g_4 g_7}{g_1 g_2}\Bigr) - \{u_1(2,4,6,1)\}_2 \otimes \Bigl(\frac{g_2 g_4}{g_1 g_7}\Bigr),
\end{equation}
with the shortcut $g_k = g_{0,k}$ for $k \in \Z/15\Z$. Let us denote by $\tilde{\xi}_{15}$ the cocycle in the right-hand side of \eqref{eq xi_15}. Lal\'\i n has shown that the cocycle $\tilde{\xi}_E$ has trivial residues \cite[Section 4.1, p.~213]{Lal15}, hence $\tilde{\xi}_{15}$ has trivial residues at the cusps.

The next task is to express $\tilde{\xi}_{15}$ in terms of the cocycles $\tilde{\xi}_1(a,b)$ with $a, b \in \Z/15\Z$. We do this using the modular complex $\mathcal{C}_{15}(3)$ from Section \ref{sec weight 3}. Using the function \verb+find_xi1ab+ from \verb+K4-modular-complex.gp+ \cite{Bru23}, we detect the following simple expression for $\tilde{\xi}_{15}$.

\begin{pro} \label{pro classes}
We have the equality of cocycles $\tilde{\xi}_{15} = - 20 \tilde{\xi}_1(1,4) + \Xi$, where $\Xi$ is a $\Q$-linear combination of coboundaries $\{u\}_2 \otimes u$ with $u \in U_2$. In particular, we have $\varphi^*(\xi_E) = -20 \xi_1(1,4)$.
\end{pro}

\section{The integration path} \label{sec path}

In Theorem \ref{thm Lalin}, the integration path $\gamma_E = \partial \Gamma$ is a closed path in $E$, and we would like to express it in terms of modular symbols on $X_1(15)$, via the modular parametrisation from Section \ref{sec parametrisation}. This is a crucial ingredient in the computation of the regulator integral on $E$. We will do this carefully in order to certify the relation (Proposition \ref{pro paths}).

Lal\'\i n \cite[Section 4.1]{Lal15} has shown that $\gamma_E$ is a generator of $H_1(E, \Z)^+$, where $(\cdot)^+$ denotes the subgroup of invariants under complex conjugation. So we first search for a generator $\gamma_{15}$ of $H_1(X_1(15),\Z)^+$. We do this with the help of SageMath \cite{SAGE}; see the notebook \verb+ModularSymbolGamma15.ipynb+ in \cite{Bru23}. For any $g \in \SL_2(\Z)$, denote by $[g] = \{g0,g\infty\}$ the associated Manin symbol, viewed in the relative homology group $H_1(X_1(15),\{\textrm{cusps}\},\Z)$. We obtain
\begin{equation} \label{eq gamma_15}
\gamma_{15} = 2 \left[\begin{pmatrix} 1 & 9 \\ 2 & 19 \end{pmatrix}\right] - \left[\begin{pmatrix} 0 & -1 \\ 1 & 11 \end{pmatrix}\right] - \left[\begin{pmatrix} 0 & -1 \\ 1 & 4 \end{pmatrix}\right] + 2 \left[\begin{pmatrix} 0 & -1 \\ 1 & 2 \end{pmatrix}\right].
\end{equation}
We therefore have $\gamma_E = \pm \varphi_*(\gamma_{15})$. The precise sign is not strictly needed in what follows, as the Mahler measure is a positive real number and the final identity will fix the sign for us. However, we want to sketch a method to determine the sign rigorously, as it could be useful in more general situations, where the integration path $\gamma$ need not be a generator of the homology group. In such a scenario, one wishes to ascertain an identity of the form $\gamma = c \cdot \varphi_*(\gamma_0)$, where $\varphi$ is the modular parametrisation, $\gamma_0$ is a modular symbol, and $c \in \Z$ is to be determined.

The idea is to integrate an invariant differential form over the cycles to be compared. By \cite[Section 4.1]{Lal15}, an invariant differential form on $E$ is given by
\begin{equation*}
\omega_E := \frac{-dx}{2(x+1)^2 (y+1) - x}.
\end{equation*}
Using \eqref{eq ut vt},  we can compute the Fourier expansion of the pull-back of $\omega_E$ to $X_1(15)$:
\begin{equation*}
W_{15}^*(\varphi^* \omega_E) = - (q-q^2-q^3+O(q^4)) \frac{dq}{q}.
\end{equation*}
A basis of $\Omega^1(X_1(15))$ is given by $\omega_{15} := 2\pi i f_{15}(\tau) d\tau$, where $f_{15} = q-q^2-q^3+O(q^4)$ is the newform of weight $2$ on $\Gamma_1(15)$. Therefore $W_{15}^*(\varphi^* \omega_E) = - \omega_{15}$. Moreover, the involution $W_{15}$ has a fixed point $\tau = i/\sqrt{15}$ in the upper half-plane, so it must act on the complex torus underlying $X_1(15)$ as $z \mapsto z_0 - z$ for some $z_0$ (it cannot be a translation). It follows that $W_{15}$ acts as $-1$ on $\Omega^1(X_1(15))$, and we conclude that $\varphi^* \omega_E = \omega_{15}$.

\begin{figure}[h]
\centering
\includegraphics[width=0.33\textwidth]{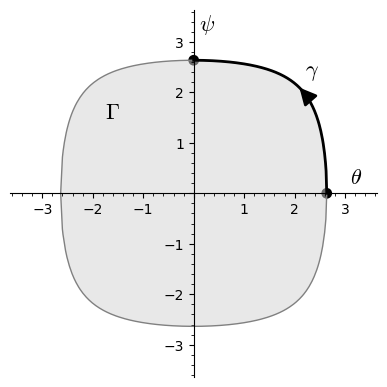}
\caption{The Deninger chain $\Gamma$, its boundary $\gamma_E$ and the path $\gamma$.}
\label{fig:gamma_E}
\end{figure}

Now let us integrate the forms $\omega_E$ and $\omega_{15}$, and compare the signs of the integrals. Following \cite[Section 4.1]{Lal15}, the path $\gamma_E$ is described using polar coordinates $x=e^{i\theta}$, $y=e^{i\psi}$ with $\theta, \psi \in [-\pi,\pi]$, and is given by the equation $\cos(\theta/2) \cos(\psi/2) = 1/4$. Since the orientation of the Deninger chain $\Gamma$ is induced by the product orientation of $[-\pi,\pi]^2$, its boundary $\gamma_E$ is oriented counterclockwise in this square (see Figure \ref{fig:gamma_E}). We can use the symmetries of $\gamma_E$ to reduce the integration path. For any automorphism $\sigma$ of $E$ defined over $\R$, we have $\sigma^* \omega_E = \varepsilon(\sigma) \omega_E$, where $\varepsilon(\sigma)=1$ if $\sigma$ preserves the orientation of $E(\R)$, and $\varepsilon(\sigma)=-1$ otherwise. Equivalently, $\varepsilon(\sigma) = 1$ if and only if $\sigma=\mathrm{id}$ or $\sigma$ has no fixed point. Applying this with the symmetries $(x,y) \mapsto (1/x,y)$ and $(x,y) \mapsto (x,1/y)$, which reverse the orientation of $E(\R)$ as well as that of $\gamma_E$, we obtain that $\int_{\gamma_E} \omega_E$ is 4 times the integral over the path $\gamma$ pictured in Figure \ref{fig:gamma_E}.

After some computation, we get
\begin{equation} \label{eq int gamma_E}
\int_{\gamma_E} \omega_E = 4 \int_\gamma \omega_E = 4 \int_0^{2 \arccos(1/4)} \frac{d\theta}{\sqrt{16 \cos^2(\theta/2)-1}} > 0.
\end{equation}
Now with the modular curve $X_1(15)$, we wish to determine the sign of $\int_{\gamma_{15}} \omega_{15}$. For this, consider the linear map $H_1(X_1(15),\{\textrm{cusps}\},\Z) \to H_1(X_1(15),\Q)$ provided by the Manin-Drinfeld theorem \cite{Dri73}. Again with SageMath, we compute that the image of $\{0,\infty\}$ is equal to $-\frac{1}{16} \gamma_{15}$ (see \verb+ModularSymbolGamma15.ipynb+ \cite{Bru23}). It follows that
\begin{equation} \label{eq int gamma_15}
\int_{\gamma_{15}} \omega_{15} = -16 \int_0^\infty \omega_{15} = 16 L(f_{15}, 1) > 0.
\end{equation}
That $L(f_{15}, 1)$ is positive can be ascertained witout much effort using the rapidly convergent series $L(f_{15}, 1) = 2 \sum_{n = 1}^{\infty} a_n e^{-2\pi n/\sqrt{15}}/n$ \cite[Proposition 7.5.8]{Coh93}. Namely, one may use the bound $|a_n| \leq n$ for $n \geq 1$, which follows from the Hasse bound on $E$ and the inspection of the coefficients $a_n$ for small $n$. Combining \eqref{eq int gamma_E} and
\begin{equation*}
\int_{\varphi_*(\gamma_{15})} \omega_E = \int_{\gamma_{15}} \varphi^* \omega_E = \int_{\gamma_{15}} \omega_{15} > 0,
\end{equation*}
we come to the following conclusion.

\begin{pro} \label{pro paths}
We have $\gamma_E = \varphi_*(\gamma_{15})$.
\end{pro}

To be fully accurate (and in order to handle more general situations), ascertaining this equality requires to compute numerically the integrals \eqref{eq int gamma_E} and \eqref{eq int gamma_15}. And since the ratio of these integrals is known to be an integer, it suffices to compute them with rigorous error bounds. The integral \eqref{eq int gamma_E} is a complete elliptic integral which can be dealt with the Arb library \cite{Joh17, Joh18}. On the other hand, \eqref{eq int gamma_15} involves integrating a modular form over a modular symbol. We can do it in the present situation thanks to the rapidly convergent series. In general, although PARI/GP \cite{PARI} and Magma \cite{Magma} can evaluate such integrals efficiently, we are not aware of implementations that prove error bounds for them.

\section{Final computation} \label{sec final}

We denote by $r_3(2)$ the Goncharov regulator map in degree $2$ for the weight $3$ complex of a curve (see \cite{Gon02}). It sends a degree $2$ cocycle to an explicit closed $1$-form on this curve. By Proposition \ref{pro paths}, we have
\begin{equation} \label{eq final 1}
\int_{\gamma_E} r_3(2)(\xi_E) = \int_{\varphi_* (\gamma_{15})} r_3(2)(\xi_E) = \int_{\gamma_{15}} r_3(2)(\varphi^* \xi_E) = \int_{\gamma_{15}} r_3(2)(\tilde{\xi}_{15}).
\end{equation}
Note that the differential form $r_3(2)(\tilde{\xi}_{15})$ is defined only on the open modular curve $Y_1(15)$. However, it has trivial residues at the cusps since the same is true for $\tilde{\xi}_{15}$, see Section \ref{sec class}. We may therefore compute the integral by choosing the representative of $\gamma_{15}$ given by \eqref{eq gamma_15}. Note that this integral involves cusps but it is absolutely convergent by \cite[Corollary 7.3]{Bru20b}. The technical details of this procedure are explained at the end of \cite[Section 8]{Bru20b}.

\begin{lem} \label{lem int coboundary}
Let $u$ be a modular unit on $X_1(N)$ such that $1-u$ is also a modular unit. For any two cusps $\alpha \neq \beta$ in $\p^1(\Q)$, we have $\int_\alpha^\beta r_3(2)(\{u\}_2 \otimes u) = \hat{\mathcal{L}}_3(u(\beta)) - \hat{\mathcal{L}}_3(u(\alpha))$, where $\hat{\mathcal{L}}_3 : \p^1(\C) \to \R$ is the single-valued trilogarithm defined in \cite[Section 2.1]{Gon02}.
\end{lem}

\begin{proof}
By \cite[Theorem 2.2]{Gon02}, we have
\begin{equation*}
r_3(2)(\{u\}_2 \otimes u) = r_3(2)(\delta(\{u\}_3)) = d r_3(1)(\{u\}_3) = d\hat{\mathcal{L}}_3(u). \qedhere
\end{equation*}
\end{proof}

Since the path $\gamma_{15}$ is closed, Lemma \ref{lem int coboundary} implies that $\int_{\gamma_{15}} r_3(2)(\{u\}_2 \otimes u) = 0$ for any $u \in U_2$. Using Proposition \ref{pro classes}, the computation \eqref{eq final 1} continues as
\begin{equation} \label{eq final 2}
\int_{\gamma_E} r_3(2)(\xi_E) = -20 \int_{\gamma_{15}} r_3(2)(\tilde{\xi}_1(1,4)).
\end{equation}
We are now in position to apply the main result of \cite{BZ23}, which computes
\begin{equation*}
\mathcal{G}(\bolda,\boldb) = \int_0^\infty r_3(2)(\tilde{\xi}(\bolda,\boldb)) \qquad (\bolda,\boldb \in (\Z/N\Z)^2),
\end{equation*}
under the assumption that the coordinates of $\bolda$, $\boldb$ and $\bolda+\boldb$ are non-zero. We may integrate along Manin symbols $[g] = \{g0,g\infty\}$ as well, noting that
\begin{equation*}
\int_{g0}^{g\infty} r_3(2)(\tilde{\xi}(\bolda,\boldb)) = \int_0^\infty r_3(2)(\tilde{\xi}(\bolda g,\boldb g)) = \mathcal{G}(\bolda g, \boldb g) \qquad (g \in \SL_2(\Z)).
\end{equation*}
Recall also that $\tilde{\xi}_1(a,b) = \tilde{\xi}((0,a),(0,b))$. Expanding \eqref{eq final 2}, we get
\begin{equation*}
\int_{\gamma_E} r_3(2)(\xi_E) = -20 \bigl(2 \mathcal{G}((2,4),(8,1)) - \mathcal{G}((1,11),(4,14)) - \mathcal{G}((1,4),(4,1)) + 2 \mathcal{G}((1,2),(4,8))\bigr).
\end{equation*}
The assumption on the coordinates of the parameters is satisfied, and \cite[Theorem 1]{BZ23} gives
\begin{equation} \label{eq final 3}
\int_{\gamma_E} r_3(2)(\xi_E) = \pi^2 L'(F,-1)
\end{equation}
with
\begin{equation} \label{eq final 4}
\begin{split}
F & = -8(G_{2,1} G_{8,-4} + G_{2,-1} G_{8,4}) + 4 (G_{1,14} G_{4,-11} + G_{1,-14} G_{4,11}) \\
& \quad + 4 (G_{1,1} G_{4,-4} + G_{1,-1} G_{4,4}) - 8(G_{1,8} G_{4,-2} + G_{1,-8} G_{4,2}).
\end{split}
\end{equation}
Here $G_{a,b}$ is a shortcut for the Eisenstein series $G^{(1);15}_{a,b}$ defined in \cite[Introduction]{BZ23} for arbitrary level $N$ by
\begin{equation*}
G^{(1);N}_{a,b}(\tau) = a_0(G^{(1);N}_{a,b}) + \sum_{\substack{m,n \geq 1 \\ (m,n) \equiv (a,b) \bmod{N}}} q^{mn/N} - \sum_{\substack{m,n \geq 1 \\ (m,n) \equiv -(a,b) \bmod{N}}} q^{mn/N} \qquad (a,b \in \Z/N\Z).
\end{equation*} 
In our situation the indices $a,b$ are non-zero modulo $15$, so that the constant terms $a_0(G_{a,b})$ vanish. The functions $G_{a,b}$ are Eisenstein series of weight $1$ on $\Gamma(15)$. Note that the products $G_\boldx G_\boldy$ appearing in \eqref{eq final 4} are actually power series in $q$, because $x_1 x_2 + y_1 y_2$ is divisible by $15$ for each such product. It follows that $F$ belongs to $M_2(\Gamma_1(15))$.

We have written a script \verb+K4-reg-Lvalue.gp+ \cite{Bru23} to automate the application of \cite[Theorem 1]{BZ23} and compute the $q$-expansion of the resulting modular form to arbitrary precision. We find that $F = -8 f_{15} + O(q^{21})$, where $f_{15}$ is the newform associated to $E$. Moreover, the Sturm bound for the space $M_2(\Gamma_1(15))$ is equal to 16 (apply \cite[Sturm's theorem, 9.4.1.2]{Ste07} with the group $\Gamma = \pm \Gamma_1(15)$, which has index $96$ in $\SL_2(\Z)$). This means that if two modular forms $F_1$ and $F_2$ in this space satisfy $F_1 = F_2 + O(q^{17})$, then $F_1 = F_2$. In our situation, this allows us to certify that $F = -8 f_{15}$. Using Theorem \ref{thm Lalin} and \eqref{eq final 3}, the Mahler measure finally equals
\begin{equation*}
m(P) = \frac{1}{4\pi^2} \int_{\gamma_E} r_3(2)(\xi_E) = \frac{1}{4\pi^2} \cdot \pi^2 L'(-8 f_{15}, -1) = -2 L'(E,-1).
\end{equation*}
This concludes the proof of Theorem \ref{thm main}.

\section*{Appendix. Tables of $3$-variable Mahler measures} \label{sec tables}

We would like to give here a list of conjectural identities for $3$-variable Mahler measures involving $L(E,3)$ for several elliptic curves $E$ over $\Q$. It is possible that our methods can be applied to prove at least some of these identities. The success of the approach will depend very much on the modular parametrisation of the elliptic curve; in our case, Proposition \ref{pro parametrisation} was crucial. This is similar to what happens for the $2$-variable Mahler measures, where the proofs using the Rogers--Zudilin method require the curve to be parametrised by modular units \cite[Section 8.4 and Chapter 9]{BZ20}.

Boyd and Rodriguez Villegas \cite{Boy05} discovered several identities of type $m(P(x,y,z)) = r \cdot L'(E,-1)$ with $r \in \Q^\times$ by looking at polynomials of the form $P = A(x)+B(x) y + C(x) z$ where $A$, $B$, $C$ are products of cyclotomic polynomials. Boyd found further examples in \cite{Boy06a, Boy06b}. We extended Boyd's search with $A$, $B$, $C$ of degree up to $5$ and found a few other examples, see Table \ref{table:1} below (we do not claim to have spotted all identities for this range of $A, B, C$). Table \ref{table:2} displays two Mahler measures which involve a combination of $L(E,3)$ and $\zeta(3)$ with polynomials $P$ of the same type. Note that $\zeta(3)$ terms also appear in \cite[Theorem 1]{BZ23}.

In the tables below, the curve $E$ is given by its Cremona label, and the integer $g$ is the genus of the Maillot variety $W_P$ (or a component of it) whose Jacobian has $E$ as an isogeny factor.

\begin{table}[h]
\centering
\begin{tabular}{c|c|c|c|c}
$P$ & $E$ & $r$ & $g$ & Source \\[5pt]
\hline
& & & & \\[-10pt]
$(x-1)^3 + (x+1) (y+z)$ & $14a4$ & $-6$ & $2$ & \cite[p.~21]{Boy06a}, \cite{Bru20a} \\[5pt]
\hline
& & & & \\[-10pt]
$(1+x)^2 + (1-x)(y+z)$ & $20a1$ & $-2$ & $1$ &  \cite{Boy06b}, \cite[p.~81]{BZ20} \\[5pt]
\hline
& & & & \\[-10pt]
$x+1 + (x-1) y + (x^2-1) z$ & $20a1$ & $-3/2$ & $3$ & \cite{Bru20a} \\[5pt]
\hline
& & & & \\[-10pt]
$(1+x)^2 (1+y) + z$ & $21a4$ & $-3/2$ & $1$ & \cite{Boy06b}, \cite[p.~81]{BZ20} \\[5pt]
\hline
& & & & \\[-10pt]
$1 + x + y - xy + z$ & $21a1$ & $-5/4$ & $1$ & \cite{Boy06b}, \cite[Section 6.3]{BZ20} \\[5pt]
\hline
& & & & \\[-10pt]
$(1+x)^2 + y + z$ & $24a4$ & $-1$ & $1$ & \cite[Section 8]{BLRVD03}, \cite{Boy06b} \\[5pt]
\hline
& & & & \\[-10pt]
$(x+1)^2 + (x^2-1) y + (x-1)^2 z$ & $48a1$ & $-2/5$ & $3$ & \cite{Bru20a} \\[5pt]
\hline
& & & & \\[-10pt]
$(x+1)^2 + (x-1)^2 y + z$ & $225c2$ & $-1/48$ & $1$ & \cite{Boy06b, Bru20a}
\end{tabular}
\caption{Conjectural identities $m(P) \stackrel{?}{=} r \cdot L'(E,-1)$}
\label{table:1}
\end{table}

\begin{table}[h]
\centering
\begin{tabular}{c|c|c|c|c|c}
$P$ & $E$ & $r$ & $s$ & $g$ & Source \\[5pt]
\hline
& & & & & \\[-10pt]
$x^2-x+1 + (x+1)(y+z)$ & $30a1$ & $-2/9$ & $-476/27$ & $1$ & \cite{Bru20a} \\[5pt]
\hline
& & & & & \\[-10pt]
$(x-1)^3 + (x+1)^3 (y+z)$ & $108a1$ & $-2/15$ & $112/15$ & $2$ & \cite[p.~21]{Boy06a}, \cite{Bru20a}
\end{tabular}
\caption{Conjectural identities $m(P) \stackrel{?}{=} r \cdot L'(E,-1) + s \cdot \zeta'(-2)$}
\label{table:2}
\end{table}

We also looked at polynomials $P(x,y,z)$ which have degree $1$ in each variable $x, y, z$, and all of whose coefficients are $\pm 1$ (or zero). It seems to be the case that every such polynomial is exact. The identities found are collected in Table \ref{table:3}. The first entry in this table is not of this shape but we include it for completeness; it already appears in \cite{BZ20}. Ringeling computed numerically the Mahler measures below, and the identities seem to hold to at least 100 digits.

A particular feature of Table \ref{table:3} is the appearance of the elliptic curve $36a1$, which has complex multiplication. The elliptic curve $450c1$ is also the first example with a curve of rank $1$.

\begin{table}[h]
\centering
\begin{tabular}{c|c|c|c|c}
$P$ & $E$ & $r$ & $g$ & Source \\[5pt]
\hline
& & & & \\[-10pt]
$(1+x)(1+y)(x+y)+z$ & $14a4$ & $-3$ & 1 & \cite[p.~81]{BZ20} \\[5pt]
\hline
&  &  &  & \\[-10pt]
$1+x+y+z+xy+xz+yz$ & $14a4$ & $-5/2$ & 1 & \cite{Bru20a} \\[5pt]
\hline
& & & & \\[-10pt]
$1+x+y+z+xy+xz+yz-xyz$ & $36a1$ & $-1/2$ & 1 & \cite{Bru20a} \\[5pt]
\hline
& & & & \\[-10pt]
$1+xy + (1+x+y)z$ & $90b1$ & $-1/20$ & 1 & \cite{Bru20a} \\[5pt]
\hline
& & & & \\[-10pt]
$(1+x)(1+y) + (1-x-y)z$ & $450c1$ & $1/288$ & 1 & \cite{Bru20a}
\end{tabular}
\caption{Conjectural identities $m(P) \stackrel{?}{=} r \cdot L'(E,-1)$}
\label{table:3}
\end{table}

\end{document}